\documentclass[a4paper, reqno]{amsart}

\allowdisplaybreaks

\usepackage[top=0.9in, bottom=0.9in]{geometry}

\usepackage[T1]{fontenc}
\usepackage[utf8]{inputenc}
\usepackage{amssymb}

\usepackage[PS]{diagrams}
\usepackage{xcolor}
\usepackage{pgf}

\newcommand{\ignore}[1]{\relax}

\newcommand{\id}{\ensuremath{\mathrm{id}}}

\newtheorem*{theorem*}{Theorem}

\theoremstyle{definition}
\newtheorem*{definition*}{Definition}

\usepackage[hypertexnames=false,colorlinks=false,pdfborderstyle={/S/U/W 1}]{hyperref}

\newcommand{\inv}{^{-1}}

\renewcommand{\rho}{\varrho}

\DeclareMathOperator*{\colim}{colim}

\newcommand{\GL}{\mathrm{GL}}

\begin{document}

\title{An elementary description of $K_1(R)$ without elementary matrices}

\date{\today}

\author{Thomas H\"uttemann}

\address{Thomas H\"uttemann\\ Queen's University Belfast\\ School of
  Mathematics and Physics\\ \break Mathematical Sciences Research Centre\\
  Belfast BT7~1NN\\ UK}

\email{t.huettemann@qub.ac.uk}

\urladdr{https://t-huettemann.github.io/}

\author{Zuhong Zhang}

\address{Zhang Zuhong\\ School of Mathematics\\ Beijing Institute Of
  Technology\\ \break 5 South Zhongguancun Street, Haidian District\\
  100081 Beijing\\ P.~R.~China}

\email{zuhong@gmail.com}


\begin{abstract}
  Let $R$ be a ring with unit. Passing to the colimit with respect to
  the standard inclusions $\GL(n,R) \rTo \GL(n+1,R)$ (which add a unit
  vector as new last row and column) yields, by definition, the stable
  linear group $\GL(R)$; the same result is obtained, up to
  isomorphism, when using the ``opposite'' inclusions (which add a
  unit vector as new first row and column). In this note it is shown
  that passing to the colimit along both these families of inclusions
  simultaneously recovers the algebraic $K$-group
  $K_1(R) = \GL(R)/E(R)$ of~$R$, giving an elementary description that
  does not involve elementary matrices explicitly.
\end{abstract}


\maketitle

Let $R$ be an associative ring with unit element~$1$, and let
$\GL(n,R)$ denote the group of invertible $n \times n$-matrices with
entries in~$R$. The usual stabilisation maps
\begin{displaymath}
  i^n_{n+1} \colon \GL(n,R) \rTo \GL(n+1,R) \ , \quad A \mapsto
  \begin{pmatrix}
    A & 0 \\
    0 & 1
  \end{pmatrix}
\end{displaymath}
are used to define the stable general linear group
$\GL(R) = \bigcup_{n \geq 3} \GL(n,R)$, or, phrased in categorical
language,
\begin{equation}
  \GL(R) = \colim \Big( \GL(3,R) \rTo^{i^3_4} \GL(4,R) \rTo^{i^4_5}
  \GL(5,R) \rTo^{i^5_6} \ldots \ \Big) \ .
  \label{eq:GL}
\end{equation}
The canonical group homomorphisms
$\iota_n \colon \GL(n,R) \rTo \GL(R)$ are injective and satisfy the
relation
\begin{equation}
  \label{eq:iota-i}
  \iota_{n+1} \circ i^n_{n+1} = \iota_n \ .
\end{equation}

There are other ``block-diagonal'' embedding
$i^{n}_j \colon \GL(n,R) \rTo \GL(n+1,R)$, for $1 \leq j \leq n+1$,
characterised by saying that the $j$th row and $j$th column of
$i^n_j (A)$ are $j$th unit vectors, and that deleting these from
$i^n_j(A)$ recovers the matrix~$A$. We will determine the result of
stabilising over first and last embeddings simultaneously, that is, we
identify the categorical colimit~$M$ of the following group-valued
infinite diagram:
\begin{equation}
  \GL(3,R) \pile {\rTo^{i^3_1} \\ \rTo_{i^3_4}}
  \GL(4,R) %
  \pile {\rTo^{i^4_1} \\ \rTo_{i^4_5}}
  \GL(5,R) %
  \pile {\rTo^{i^5_1} \\ \rTo_{i^5_6}}
  \cdots %
  \pile {\rTo^{i^{n-1}_1} \\ \rTo_{i^{n-1}_n}} %
  \GL(n,R) %
  \pile {\rTo^{i^n_1} \\ \rTo_{i^n_{n+1}}} %
  \GL(n+1,R) %
  \pile {\rTo^{i^{n+1}_1} \\ \rTo_{i^{n+1}_{n+2}}} %
  \cdots \ .
  \label{eq:diag}
\end{equation}
By the general theory of colimits, the group $M$ comes equipped with
canonical group homomorphisms $\alpha_n \colon \GL(n,R) \rTo M$
satisfying the relations
\begin{equation}
  \label{eq:relation_alpha}
  \alpha_{n+1} \circ i^n_j = \alpha_n \quad \text{(\(j = 1,\, n+1\))}
  \ . 
\end{equation}

\begin{theorem*}
  The group~$M$ is canonically isomorphic to~$K_1(R)$.
\end{theorem*}

\begin{proof}
  First we observe that in~$M$ we have the commutation relation
  \begin{equation}
    \label{eq:commute}
    \alpha_n(X) \alpha_n(Y) = \alpha_n(Y) \alpha_n(X) \quad \text{for
      all } X, Y \in \GL(n,R) \ .
  \end{equation}
  Indeed, by~\eqref{eq:relation_alpha} we can re-write
  \begin{displaymath}
    \alpha_n(X) = \alpha_{2n} \big( i^{2n-1}_{2n} i^{2n-2}_{2n-1} \ldots
    i^n_{n+1} (X) \big) %
    \quad \text{and} \quad %
    \alpha_n(Y) = \alpha_{2n} \big( i^{2n-1}_{1} i^{2n-2}_{1} \ldots
    i^n_{1} (Y) \big) \ ,
  \end{displaymath}
  and the arguments of $\alpha_{2n}$ are block-diagonal matrices
  of the form
  \begin{displaymath}
    \begin{pmatrix}
      X & 0   \\
      0 & I_n \\
    \end{pmatrix}
    \quad \text{and} \quad
    \begin{pmatrix}
      I_n & 0 \\
      0 & Y   \\
    \end{pmatrix}
  \end{displaymath}
  which commute in~$\GL(2n,R)$; hence their images under~$\alpha_{2n}$
  must commute as well.

  \medbreak

  Let $E(n,R)$ denote the subgroup of~$\GL(n,R)$ generated by the
  elementary matrices \cite[\S{}V.1]{MR0249491}. Since
  $E(n,R) = \big[ E(n,R),\ E(n,R) \big]$ for all $n \geq 3$,
  cf.~\cite[Corollary~V.1.5]{MR0249491}, the commutation
  relation~\eqref{eq:commute} implies
  \begin{equation}
    E(n,R) \subseteq \ker (\alpha_n) \ . \label{eq:En_ker} 
  \end{equation}

  \medbreak

  Since the diagram~\eqref{eq:GL} defining $\GL(R)$ is contained in
  the diagram~\eqref{eq:diag} defining~$M$ there is a canonical group
  homomorphism $\alpha \colon \GL(R) \rTo M$ described completely by
  $\alpha \circ \iota_n = \alpha_n$, that is, $\alpha|_{\GL(n,R)} =
  \alpha_n$.

  \medbreak

  Let $E(R) = \bigcup_{n \geq 3} E(n,R) \subseteq \GL(R)$ be the
  stabilisation \textit{via} the embeddings~$i^n_{n+1}$. In view
  of~\eqref{eq:En_ker} above we have the inclusion
  \begin{equation}
    E(R) = \bigcup_{n \geq 3} E(n,R) \subseteq \ker \alpha %
    \ .
    \label{eq:ker}
  \end{equation}

  \medbreak

  The group $E(R)$ is normal in~$\GL(R)$,
  cf.~\cite[Theorem~V.2.1]{MR0249491}, and $K_1(R) = \GL(R)/E(R)$ is
  an \textsc{abel}ian group \cite[p.~229]{MR0249491}. We write
  $\pi \colon \GL(R) \rTo K_1(R)$ for the canonical projection. Let
  $\pi_n = \pi \circ \iota_n$ denote the restriction of~$\pi$
  to~$\GL(n,R)$, and write $[X] = \pi_n(X)$ for the class of
  $X \in \GL(n,R)$ in~$K_1(R)$. By~\eqref{eq:ker} we obtain a
  factorisation $\lambda \colon K_1(R) \rTo M$ of~$\alpha$ with
  $\lambda \circ \pi = \alpha$. Explicitly, $\lambda$ is described by
  the formula
  \begin{equation}
    \lambda \colon [X] = \pi_n(X) \ \mapsto \ \alpha \circ \iota_n (X)
    = \alpha_n(X) \ , \quad \text{for } X \in \GL(n,R) \ .
    \label{eq:lambda_explicit}
  \end{equation}

  \medbreak

  We observe the relation $\pi_{n+1} \circ i^n_{j} = \pi_n$.  Indeed,
  for $X \in \GL(n,R)$ the matrices $i^n_{j} (X)$ and $i^n_{n+1}(X)$
  are related by the expression
  \begin{displaymath}
    i^n_{j} (X) = P\inv i^n_{n+1}(X) P
  \end{displaymath}
  for a permutation matrix $P \in \GL(n+1,R)$. It follows that said
  two matrices have the same image under~$\pi_{n+1}$ in the
  \textsc{abel}ian group~$K_1(R)$ whence, using~\eqref{eq:iota-i},
  \begin{displaymath}
    \pi_{n+1} \circ i^n_j (X) = \pi_{n+1} \circ i^n_{n+1} (X) = \pi
    \circ \iota_{n+1} \circ i^n_{n+1} (X) = \pi \circ \iota_n (X) =
    \pi_n(X) \ .
  \end{displaymath}
  The various maps~$\pi_n$ thus form a ``cone'' on the
  diagram~\eqref{eq:diag} and induce a map $\rho \colon M \rTo K_1(R)$
  such that
  \begin{equation}
    \label{eq:beta_n_is_rho_alpha_n}
    \pi_n = \rho \circ \alpha_n \ .
  \end{equation}

  \medbreak

  We verify the equality $\lambda \circ \rho = \id_M$. By the
  universal property of colimits, it is enough to show that
  $\lambda \circ \rho \circ \alpha_n = \alpha_n$ for all~$n$. But for
  $X \in \GL(n,R)$ we calculate
  \begin{displaymath}
    \lambda \circ \rho \circ \alpha_n (X)
    \underset{\eqref{eq:beta_n_is_rho_alpha_n}}= \lambda \circ \pi_n
    (X) \underset{\eqref{eq:lambda_explicit}}= \alpha_n(X) \ ,
  \end{displaymath}
  using relation~\eqref{eq:beta_n_is_rho_alpha_n} and the explicit
  description~\eqref{eq:lambda_explicit} of~$\lambda$ above.

  Using the same relations again, in opposite order, we finally verify
  that $\rho \circ \lambda = \id_{K_1(R)}$. Let $X \in \GL(n,R)$
  represent the element $[X] \in K_1(R)$ as before; then
  $\rho \circ \lambda \big( [X] \big)
  \underset{\eqref{eq:lambda_explicit}}= \rho \circ \alpha_n (X)
  \underset{\eqref{eq:beta_n_is_rho_alpha_n}}= \pi_n (X) = [X]$.
\end{proof}

With minor changes the argument also shows that the the colimit of the
group-valued diagram
\begin{equation*}
  \GL(3,R) \pile {\rTo^{i^3_1} \\ \rTo \\ \rTo \\ \rTo_{i^3_4}}
  \GL(4,R) %
  \pile {\rTo^{i^4_1} \\ \rTo \\ \rTo \\ \rTo \\ \rTo_{i^4_5}}
  \GL(5,R) %
  \pile {\rTo^{i^5_1} \\ \rTo \\ \rTo \\ \rTo \\ \rTo \\ \rTo_{i^5_6}}
  \cdots %
  \pile {\rTo^{i^{n-1}_1} \\ \vdots \\ \ \\ \rTo_{i^{n-1}_n}} %
  \GL(n,R) %
  \pile {\rTo^{i^n_1} \\ \vdots \\ \ \\ \rTo_{i^n_{n+1}}} %
  \GL(n+1,R) %
  \pile {\rTo^{i^{n+1}_1} \\ \vdots \\ \ \\ \rTo_{i^{n+1}_{n+2}}} %
  \cdots
\end{equation*}
is canonically isomorphic to~$K_1(R)$.


\begin{thebibliography}{Bas68}

\bibitem[Bas68]{MR0249491}
Hyman Bass.
\newblock {\em Algebraic {$K$}-theory}.
\newblock W. A. Benjamin, Inc., New York-Amsterdam, 1968.

\end{thebibliography}
\end{document}